\documentclass[a4paper]{amsart}
\usepackage[utf8]{inputenc}
\usepackage{hyperref}
\usepackage{amsmath,amsfonts,amssymb}
\usepackage{systeme}
\usepackage{bbm}
\usepackage{listings}
\usepackage{graphicx}
\usepackage{natbib}
\usepackage{bm}
\mathchardef\mhyphen="2D
\usepackage{float}
\usepackage{color}
\usepackage[dvipsnames]{xcolor} 
\usepackage{commath}
\usepackage{amssymb}
\usepackage{bm}
\usepackage{ mathrsfs }
\usepackage[colorinlistoftodos,prependcaption,textsize=tiny]{todonotes} 
\usepackage{ulem} 

\definecolor{webgreen}{rgb}{0,.5,0} 

\definecolor{darkblue}{rgb}{0,0,0.3}
\definecolor{lightblue}{rgb}{0.9,0.9,1}
\definecolor{blueKP}{HTML}{3F83DF}

\newtheorem{theorem}{Theorem}
\newtheorem{proposition}{Proposition}
\newtheorem{lemma}{Lemma}
\newtheorem{corollary}{Corollary}

\newtheorem{definition}{Definition}

\definecolor{light-gray}{gray}{0.95}  
\usepackage{marginnote} 
\title[The switch process and geometric divisibility]{ Characteristics of the switch process and geometric divisibility}
\usepackage{amsaddr} 
\author[H.Bengtsson]{Henrik Bengtsson }
\address{Department of Statistics,\\
Lund University}
\email{Henrik.Bengsson@stat.lu.se, \href{https://orcid.org/0000-0002-9280-4430}{ \includegraphics[height=2.2mm]{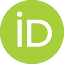} 0000-0002-9280-4430}}  
\date{4 May, 2022}

\bgroup

\begin{document}
\begin{abstract} 
The switch process alternates independently between $1$ and $-1$, with the first switch to $1$ occurring at the origin. The expected value function of this process is defined uniquely by the distribution of switching times. The relation between the two is implicitly described through the Laplace transform, which is difficult to use for determining if a given function is the expected value function of some switch process. We derive an explicit relation under the assumption of monotonicity of the expected value function. It is shown that geometric divisible switching time distributions correspond to a non-negative decreasing expected value function.
Moreover, an explicit relation between the switching time distribution and the autocovariance function of the switch process stationary counterpart is obtained, which allows parallel results for the stationary switch process and its covariance function. These results are applicable to approximation methods in statistical physics, and an example is presented. 
\end{abstract} 

\keywords{switch process, renewal theory, geometric divisibility, binary processes} 

\maketitle

\section{Introduction} \noindent
The study of binary stochastic processes has a long-standing tradition in probability theory. There exist many versions of such processes, for example, the telegraph process in continuous time or the simple discrete-time Markov chain. These processes found applications in many fields, for example, in renewal theory, signal processing, \citep{SBRS}, and in statistical physics, \citep{IIA}. 

The focus of this paper is the switch process with independent switching times. More specifically, we consider a continuous-time stochastic process taking values in $\{-1,1\}$, starting at $1$ at the origin, and then switching according to an i.i.d sequence of non-negative random variables. The switch process always starts from one hence it is not stationary, however, a convenient stationary counterpart can be defined. This counterpart will be referred to as the stationary switch process. 

The expected value of the switch process is intrinsically connected with the switching time distribution. This is also the case for the covariance of the stationary switch process. Formalizing this connection is the main contribution of the paper, among other contributions such as formulation and deriving the underlying properties of the switch process. The connection also leads to a class of distributions which constitutes a proper sub-classes of geometric infinite divisible distributions introduced by \cite{klebanov85}. 

There is a natural connection between geometric divisibility of the switching time distribution and thinning of a renewal process. The main result of this paper allows also to recover the original renewal arrival times distribution from the thinned process, see also \cite{Thin}. 

Obtaining the switching time distribution from the covariance function of a stationary binary process is a well-known problem in signal processing \citep{SBRS}. In this work, a complete solution is presented for the class of geometric divisible distributions. This is done using a relation between the covariance of the stationary switch process and the expected value of the switch process. Therefore an analogous characterization is provided to explicitly identify a function as the expected value of the switch process. 

Approximating exceedance time distributions is important in statistical physics because the explicit analytical solution is a well-known difficult problem in probability theory. There are several possible approaches and among them are numerical based on the generalized Rice formula, see for example \cite{Joint22}, and quasi-analytical such as the independent interval approximation (IIA) framework, see for example \cite{Sire2007}. Our derived relation between expected value, covariance and the switching time distribution contributes to this framework. 

Often, the tail behavior is of the main interest, see for example \cite{IIA} and \cite{Sire2007}. There is a challenge from a mathematical point of view to provide conditions when this approximation leads to a valid probability distribution, which is not obvious at all. The results of this paper allow us to obtain an explicit representation of the approximated distribution for many stochastic processes commonly used in statistical physics. This does not only provides information about the tail behavior but also yields the explicit approximated distribution of excursions above or below zero. The diffusion process in two dimensions is used to illustrate this application to IIA approximation of the exceedance time distribution. 

The paper is self-contained in that it contains the formulation of the switch process, complete derivations of the relations between the expected value and covariance of switch process, and the switching time distribution. 
The structure of the paper is as follows. In Section~\ref{preliminaries}, preliminaries such as the definition of the switch process, results connected with its expected value, and the concept of $r$-geometric divisibility are presented. The relation between the $2$-geometric divisible switching time distributions and the expected value of the switch process is presented in Section~\ref{mainresult}. The definition of the stationary switch process and the relation between its covariance and the expected value function of the switch process are presented in Section~\ref{SBRS} and are applied to obtain a parallel characterization of the covariance of the switch process. Some examples of deriving expected value and covariance for different switching time distributions follow in Section~\ref{examples}.  In Section~\ref{appli}, an application to the independent interval approximation (IIA) method in statistical physics is presented. 

\section{Preliminaries} \label{preliminaries}
\subsection{The switch process and its expected value}
A natural definition of the switch process is through the counting process. Let $\left\{ T_k \right\}_{k\geqslant 1}$ be a sequence of i.i.d non-negative random variables with the distribution function $F$. Additionally, $F$ is assumed to be absolutely continuous with respect to the Lebesgue measure with the density that is bounded on any closed interval of the positive half-line, and such that $F(0)=0$. Define a count process for $t\in [0,\infty)$, by 

\begin{align*}
N(t)=
\begin{cases}
\sup \left\{ n\in \mathbb N;  { \sum_{k=1}^n {T_k}}  \leqslant t   \right\}, & t\geqslant T_1,\\ 
0, & 0\le t<T_1,
\end{cases}
\end{align*}
In other words, $N(t)$ is the number of renewal events up to a time point t. 
 
\begin{definition} \label{Defswitch}
Let $N(t), t\ge 0$ be a count process. Then the switch process is defined by 
\begin{align*}
X(t)=(-1)^{N(t)},~t\ge 0.
\end{align*}
\end{definition} \noindent
The process $X(t)$ switches between the values $1$ and $-1$ at each renewal event, hence the name. 

Before moving to the expected value function of $X(t)$, a comment concerning convolutions is needed. Throughout this paper, $\star$ will denote convolutions of probability measures, i.e. the distribution functions of the sum of i.i.d random variables. Convolutions of functions will be denoted by $\ast$. The switch process and all its properties are determined by the switching time distribution $F$. One such property is its expected value function, that throughout this paper will be denoted by $E(t)=EX(t)$.
The following result is well-known in renewal theory, however, for completeness, a proof is provided in Appendix~\ref{proofs}, together with the other proofs relating to this section. 
\begin{proposition} \label{EXt}
Let $X(t)$ be a switch process. Then its expected value function is given through the switching time distribution $F$ by 
\begin{align*} 
E(t)&=\sum_{n=0}^\infty (-1)^n \left( F^{n \star } - F^{(n+1) \star } \right)(t). 
\end{align*}
\end{proposition}
\noindent  
The implicit nature of this relation between the expected value function and the switching time distribution poses a problem when one wants to verify that a given function is the expected value function of a switch process. To deal with the repeated convolutions we use the Laplace transform, which is one of the standard tools in renewal theory. In the Laplace domain, the relation between the switching time distribution and the expected value function becomes explicit. To clarify the notation of the next proposition, let $\Psi_F(s)=\int_0^\infty e^{-ts}dF(t)$, i.e $\Psi_F$ is the Laplace transform of the distribution function $F$. Throughout the paper $\mathcal{L}(\cdot)$ is used to denote Laplace transform. 

\begin{proposition} \label{LEXt}
Let $X(t)$ be a switch process with the expected value function $E(t)$. Then for $s>0$,
\begin{align*}
\mathcal{L}(E)(s)&=\frac{1}{s}\frac{1-\Psi_F(s)}{1+\Psi_F(s)}, \\
\Psi_F(s)&=\frac{1-s\mathcal{L}(E)(s)}{1+s\mathcal{L}(E)(s)}. 
\end{align*}
\end{proposition}
\noindent
This result can initially appear to be explicit when compared to Proposition~\ref{EXt}, but still hard to utilize in practice. The right-hand side of the last expression in Proposition~\ref{LEXt} must be a completely monotone function. Recall that a function is completely monotone if and only if $(-1)^n \frac{d^n}{ds^n} f(s)\geqslant0$ for $s>0$ and all $n\geqslant 0$. This is a consequence of Bernstein’s theorem (see Theorem 1, page 415 in \cite{FellerV2}) which states that a function is the Laplace transform of a probability distribution if and only if it is completely monotone. We conclude that a simple criterion that characterizes the set of functions such that they correspond to expected value functions of some switch process is not easy to obtain from its definition.

The limiting properties of the expected value of the switch process are important. The effect of starting from 1 at time zero will diminish and for large $t$ the expected value tends to zero due to the symmetry, i.e. equal chances to take value one and minus one. This leads to the next proposition which follows from the Key Renewal Theorem, see Appendix~\ref{proofs}.

\begin{proposition} \label{limitXt}
Let $X(t)$ be a switch process, with a continuous switching time distribution and the expected value $E(t), ~t\ge 0$, then 
\begin{align*}
\lim_{t \rightarrow 0^+} E(t)=1, \ \ \lim_{t \rightarrow \infty} E(t)=0.
\end{align*}
\end{proposition}
\noindent
In characterizing the functional properties of the expected value function of the switch process, the derivative of $E(t)$ with respect to time is needed.

\begin{proposition} \label{EXtprim}
Let $X(t)$ be a switch process with the expected value function $E(t)$. Let the switching time distribution have a density, $f$, with support in $[0,\infty)$ and that there exists $l\in \mathbb{N}$ for which $\sup_{t>0}f^{\ast l}(t)<\infty$. Then 
\begin{align*}
E'(t)=2\sum_{k=1}^{\infty} (-1)^k f^{\ast k}(t), t>0. 
\end{align*}
Moreover, the convergence of the above series as well as the infinite series in Proposition~\ref{EXt}, is locally uniform.
\end{proposition} \noindent
With the expected value functions and its properties presented, we can now present some properties of a special class of switching time distributions. 

\subsection{Geometric divisibility} \noindent 
\cite{klebanov85} introduced the concept of geometric infinite divisibility. It describes distributions that can be represented as a sum of i.i.d random variables where the number of terms in the sum follows a geometric distribution with an arbitrary parameter $p\in (0,1)$. For a broader overview of geometric infinite divisibility see \cite{TomonGID}, and \cite{onGID}. Our main focus is on a weaker concept than the geometric infinite divisibility, defined next. 

\begin{definition}\label{rdiv}
Let $\nu_p$ be a geometric random variable with the probability mass function $p_{\nu_p} (k)=(1-p)^{k-1}p$ for $k=1,2...$ and  $\{\tilde{W}_k\}_{k\geqslant 1}$ a sequence of i.i.d non-negative random variables independent of $\nu_p$. If a random variable has the stochastic representation 
\begin{align*}
W=\sum_{k=1}^{\nu_p} \tilde{W}_k,
\end{align*}
then $W$ follows a $r$-geometric divisible distribution with $r=E\nu_p$ and said to belong to the class $GD(r)$ and we write $F\in GD(r)$
\end{definition}
\noindent
The distribution of $\tilde{W}$ is then called the $r$-geometric divisor the distribution of $W$. The notion of divisibility comes naturally from splitting a random variable into a sum of smaller parts, in essence dividing it into, on average, $r$ random variables. If $W$ belongs to $GD(r)$ with the divisor $\tilde{W}$, it follows from Wald's equation that 
\begin{align*}
EW=rE\tilde{W}.  
\end{align*} 
\noindent
If $X$ belongs to $GD(r)$, with the distribution function $F$, then it has the following Laplace transform
\begin{align*}
\Psi_F(s)&=\frac{\frac{1}{r}\Psi_{\tilde{F}}(s)}{1-(1-\frac{1}{r})\Psi_{\tilde{F}}(s)},
\end{align*}
where $\Psi_{\tilde{F}}$ is the Laplace transform of the $r$-geometric divisor $\tilde{W}$. 
By solving the expression above for $\Psi_{\tilde{F}}$, conditions are obtained to verify if a given distribution is $r$-geometric divisible. The function obtained should be a completely monotone function which again follows from Bernstein’s theorem. 
To highlight this we have the next proposition. 

\begin{proposition} \label{LGD} A distribution $F$ belongs to $GD(r)$, $r>1$, if and only if
\begin{align*}
\frac{r\Psi_F(s)}{1 +(r-1)\Psi_F(s)}
\end{align*} 
is a completely monotone function, for $s \geqslant 0$ and is equal to one for $s=0$. 
\end{proposition}

For the class of geometric infinite divisible distributions, the divisibility property needs to be satisfied for all $r>1$. The concept of $r$-geometric divisibility is hence weaker. 
Throughout the rest of the paper, $GD(\infty)$ will refer to the set of geometric infinite divisible distributions, and $GD(r)$ denotes the set of $r$-geometric divisible distributions. There exist a monotone relation between sets $r$-geometric divisible distributions. 

\begin{proposition} \label{setlemma}
For $1<u\leqslant r < \infty$, then $GD(r) \subseteq GD(u)$. 
\end{proposition}
\noindent
In essence if $W$ belongs to $GD(u)$ it also belong to $GD(r)$, for $u\leqslant r$. 

\section{Switch processes with monotonic expected value function} \label{mainresult} \noindent
In this section, we fully characterize switch processes with monotonic expected value functions. For the main result, we recall the assumptions on the switching time distribution. First, it is assumed here that $F(t)$ has support on $(0,\infty)$ and that its corresponding density $f(t)$ is well defined at each $t\in(0,\infty)$. Secondly, it is assumed that $F(0)=0$, i.e, there is zero probability of instantaneous switching back after a switch.     
\begin{theorem} \label{Th1}
Let $X(t)$ be a switch process, with $E(t)$ as its expected value function and with the switching time distribution $F(t)$. Then the following conditions are equivalent 
\begin{itemize}
    \item[$(i)$] $F(t)\in GD(2)$
    \item[$(ii)$] $E(t)$ is non-negative and decreasing. 
\end{itemize}
\end{theorem}

\begin{proof}
\textit{$(i) \Rightarrow (ii)$:} \\
Since $F(t)\in GD(2)$, it has the following Laplace transform, as described in Section~\ref{preliminaries},
\begin{align*}
    \Psi_F(s)=\frac{\frac{1}{2} \Psi_{\tilde{F}(s)}}{1- \frac{1}{2} \Psi_{\tilde{F}(s)}}.
\end{align*}
Substituting this into the first expression of Proposition~\ref{LEXt}, we have
\begin{align*}
    \mathcal{L}(E)(s)&=\frac{1}{s}\frac{1-\frac{\frac{1}{2} \Psi_{\tilde{F}}(s)}{1- \frac{1}{2} \Psi_{\tilde{F}}(s)}}{1+\frac{\frac{1}{2} \Psi_{\tilde{F}}(s)}{1- \frac{1}{2} \Psi_{\tilde{F}}(s)}} 
    =\frac{1}{s} (1- \Psi_{\tilde{F}}(s)), 
\end{align*}
which is equivalent to  
\begin{align*}
    s\mathcal{L}(E)(s)-1&=-\Psi_{\tilde{F}}(s).
\end{align*}
Using the Laplace transform of derivative, $\mathcal{L}(h')(s)=s\mathcal{L}(h)(s)-h(0)$, and Proposition~\ref{limitXt},
\begin{align*}
     \mathcal{L}(-E')(s)&=\Psi_{\tilde{F}}(s).
\end{align*}
By taking the inverse Laplace transform, this implies that $-E'(t)$ is a probability density function. Therefore, to satisfy the limiting results of Proposition~\ref{limitXt}, $E(t)$ must satisfy the conditions of $(ii)$. \\

\noindent
\textit{$(ii) \Rightarrow (i)$:} \\ 
Under the assumptions of $(ii)$ and the results of Proposition~\ref{limitXt} we have
\begin{align*}
    \int_0^\infty E'(t)dt=\lim_{t \rightarrow \infty} E(t) - \lim_{t \rightarrow 0} E(t)=-1,
\end{align*}
$-E'(t)$ is thus a probability density function. Combining the derivative property of Laplace transform, presented above, and Proposition~\ref{limitXt} into the second equation of Proposition~\ref{LEXt}, we have 
\begin{align*}
    \Psi_F (s)&=\frac{1-s\mathcal{L}(E)(s)}{1+s\mathcal{L}(E)(s)} 
    =\frac{\mathcal{L}(-E')(s)}{2-\mathcal{L}(-E')(s)} 
    =\frac{\frac{1}{2} \mathcal{L}(-E')(s) }{1-\frac{1}{2} \mathcal{L}(-E')(s)}.
\end{align*}
This is the Laplace transform of a $GD(2)$ distribution, as described in Section~\ref{preliminaries}. Therefore $F(t)\in GD(2)$, with the divisor $-E'(t)$ which yields $(i)$.
\end{proof}

Theorem~\ref{Th1} directly relates functional properties of the expected value of the switch process with the switching time distribution for the class of $GD(2)$ distributions. In many situations, there is a need to obtain either the expected value function from the switching time distribution or the switching time distribution from the expected value function. The first problem is directly solved by using Proposition~\ref{LEXt}, but for the second problem, it is not usually known what kind of expected value function will produce proper switching time distributions, if any at all. By combining Theorem~\ref{Th1} and properties of $E(t)$ derived in Section~\ref{preliminaries}, a partial solution can be obtained for the case when switching time distribution belongs to $GD(2)$. To highlight this partial characterization we have the following proposition which follows directly from the second part of the proof of Theorem~\ref{Th1}.  
 
\begin{proposition} \label{Efunk}
Let $E(t)$ be a function for $t\geqslant0$ such that the following conditions are satisfied
\begin{itemize}
    \item[$(i)$] $\lim_{t\rightarrow 0+}E(t)=1$,
    \item[$(ii)$] $\lim_{t\rightarrow \infty}E(t)=0$,
    \item[$(iii)$] $E(t)$ is at least once differentiable on $(0,\infty)$, 
    \item[$(iv)$] $E'(t)\leqslant 0$, for all $t\geqslant0$,
\end{itemize}
then it is an expected value function of a switch process with a $GD(2)$ switching time distribution.
\end{proposition} 

From Theorem~\ref{Th1} an explicit representation of the $2$-geometric divisors distribution function is obtained. 

\begin{corollary} \label{Corr1}
Let the switching time distribution, $F(t)$, belong to $GD(2)$, with the divisor $\tilde{F}(t)$, then for $t\geqslant 0$
\begin{align*}
    E(t)&=1-\tilde{F}(t), \\
    E'(t)&=-\tilde{f}(t), 
\end{align*}
\end{corollary}

Corollary~\ref{Corr1} gives an explicit representation of the distribution function and density for the $2$-geometric divisor, of the switching time distribution in terms of $E(t)$. The standard method to obtain the switching time distribution from the divisor is by using the Laplace transform as seen in Section~\ref{preliminaries}. 

Proposition~\ref{setlemma} can be used to extend the results of Theorem~\ref{Th1}. 
\begin{corollary}\label{Corr2} 
Let the switching time distribution be $GD(r)$, for some $r\geqslant 2$, then the corresponding expected value function of the switch process, $E(t)$, is non-negative and decreasing for $t\geqslant 0$. 
\end{corollary}
\noindent
However, a non-negative and decreasing expected value function does not, necessarily imply a $r$-geometric divisible, switching time for $r>2$.

Let us consider a switch process  which is constructed from a count process $N(t)$ and satisfying the conditions of Theorem~\ref{Th1}. Further, let $\tilde{N}(t)$ be a count process with the arrivals times distributed according to the divisor of this switch process. The two count processes are related through thinning. More specifically, $N(t)$ is a thinning of $\tilde{N}(t)$, with the probability of thinning equal to $\frac{1}{2}$. 
Thus we have the following result.
\begin{corollary}
A switch process $X(t)$ is $\frac 1 2 $-thinned if and only if its expected value is non-negative and decreasing. 
\end{corollary}
From a given trajectory of $N(t)$, the trajectory of process $\tilde{N}(t)$ cannot be recovered, in general. However, it follows from Proposition~\ref{Efunk} that the distribution of arrival times of $\tilde{N}(t)$ can be recovered. For further relations between geometric divisibility of the switching time distribution  and the thinned renewal processes, see \cite{Thin,ThinG}. 
  
\section{The autocovariance of the stationary switch process} \label{SBRS} \noindent
The switch process is not stationary, however there exist a stationary counterpart. In order to define it, the behavior around zero needs to be addressed. Let $\mu$ be the expected value of the switching time distribution, and $((A,B),\delta )$ be non-negative random variables, mutually independent and independent of $X(t)$ with the following densities
\begin{align*}
    f_{A,B}(a,b)&=\frac{f_T(a+b)}{\mu}, \\
    f_A(t)&=f_B(t)=\frac{1-F_T(t)}{\mu}, \\
    P\{\delta=1\}&=P\{\delta=-1\}=\frac{1}{2}.
\end{align*}
\begin{definition} \label{Defssp}
Let $X_+(t)$ and $X_-(t)$ be two independent switch processes, and $((A,B),\delta )$ be as described above, define
\begin{align*}
Y(t)=
\begin{cases}
- \delta, & -B < t <A, \\
\delta X_+(t-A), &  t \geqslant A, \\
-\delta X_-(-(t+B)), & t \leqslant -B.  
\end{cases}
\end{align*}
Then $Y(t)$ called a stationary switch process. 
\end{definition}
\noindent
The stationarity of $Y(t)$ follows from well known results in renewal theory, which follows from the key renewal theorem, see for example \cite{theoryofpoint2}.

The switch process is characterized by its expected value function $E(t)$ or, equivalently, by its switching time distribution. Obviously, the stationary switch process is also characterized by the switching time distribution. It will be seen that it is also characterized by its covariance function, denoted by $C(\cdot)$. This is due to the following relation between the expected value function and the covariance function.

\begin{proposition} \label{propcovexp}
Let $Y(t)$ be a stationary switch process, $E(t)$ be the expected value of the corresponding switch process, and let $\mu$ be the expected value of the switching time distribution. Then for $s>0$  
\begin{align*} 
    \mathcal{L}(C)(s)=\frac{1}{s} \left( 1-\frac{2}{\mu }\mathcal{L}(E)(s) \right).
\end{align*}
\end{proposition}

\begin{proof}
Starting with the covariance of $Y(t)$, and utilizing symmetry, we have \\ $(-Y(t)\vert \delta=1) \stackrel{d}{=}(Y(t)\vert \delta=-1)$, and for $t>0$
\begin{align*}
    C(t)&=E \left( E(Y(t)Y(0) \vert \delta)\right) \\
    &=\frac{1}{2}\left( E(Y(t) (-\delta) \vert \delta=1) + E(Y(t) (-\delta) \vert \delta=-1 )\right) \\
    &=-E(Y(t) \vert \delta=1) \\    
    &=-\left( \int_0^\infty E(Y(t) \vert \delta=1, A=x) f_{A \vert \delta=1}(x)dx \right)\\
    &= - \left( \int_0^t E(\delta X(t-x)\vert \delta=1, A=x)f_A(x)dx + \int_t^\infty (-1) f_A(x)dx \right) \\
    &= -\int_0^t E(t-x)f_A(x)dx + 1-F_A(t). 
\end{align*}
Since $E(t-x)=0$, for $x>t$, we obtain
\begin{align*}
    C(t)=1-F_A(t)-\left(E\ast f_A \right)(t). 
\end{align*}
\noindent
From Definition~\ref{Defssp} and the Laplace transform of $F_A(t)$ 
\begin{align*}
    \mathcal{L}(f_A)(s)=\frac{1-\Psi_F(t)}{s \mu}.
\end{align*}
Using the above expression and Proposition~\ref{LEXt}, 
\begin{align*}
     \mathcal{L}(C)(s)&= \frac{1}{s} - \frac{1}{s}\mathcal{L}(f_A)(s)-\mathcal{L}(E)(s) \mathcal{L}(f_A)(s) \\ 
    &=\frac{1}{s}-\mathcal{L}(f_A)(s) \left(\frac{1}{s}+\frac{1}{s}\frac{1-\Psi_F(s)}{1+\Psi_F(s)}  \right)\\
    &=\frac{1}{s} - \mathcal{L}(f_A)(s) \left( \frac{2}{s}\frac{1}{1+\Psi_F(s)} \right) \\
    &=\frac{1}{s} - \frac{1-\Psi_F(s)}{\mu s} \left( \frac{2}{s}\frac{1}{1+\Psi_F(s)} \right) \\
    &=\frac{1}{s} - \frac{2}{\mu s} \left( \frac{1}{s} \frac{1-\Psi_F(s)}{1+\Psi_F(s)} \right)\\
    &=\frac{1}{s}-\frac{2}{\mu s}\mathcal{L}(E)(s)=\frac{1}{s}\left( 1- \frac{2}{\mu }\mathcal{L}(E)(s) \right).
\end{align*}
\end{proof}
The above relation between $E(t)$ and $C(t)$ is somewhat implicit. By investigating the limit when $t$ approaches zero of $C(t)=1-F_A(t)-\left(E\ast f_A \right)(t)$, it is clear that $C(0)=1$. With this remark, the relation between $E(t)$ and $C(t)$ becomes explicit in the next theorem. 

\begin{theorem} \label{Th2}
Let $C(t)$ be the covariance of the stationary switch process, $E(t)$ be the expected value function of the switch process, and $\mu$ be the expected value of the switching time distribution, then for $t\geqslant0$ 
\begin{align*}
    C'(t)=-\frac{2}{\mu}E(t).
\end{align*}
\end{theorem}

\newpage
\begin{proof}
From Proposition~\ref{propcovexp} we have 
\begin{align*}
    \mathcal{L}(C)(s)&=\frac{1}{s}\left(1-\frac{2}{\mu }\mathcal{L}(E)(s)\right).
\end{align*}
Using the the property of the Laplace transform that $\mathcal{L}(f')=s\mathcal{L}(f)-f(0)$ and $C(0)=1$,
\begin{align*}
    \mathcal{L}(C)(s)&=\frac{1}{s}-\frac{2}{\mu s}\mathcal{L}(E)(s), \\
    s\mathcal{L}(C)(s)-1&=-\frac{2}{\mu}\mathcal{L}(E)(s), \\
    \mathcal{L}(C')(s)&=-\frac{2}{\mu}\mathcal{L}(E)(s),\\
    C'(t)&=-\frac{2}{\mu}E(t).
\end{align*}
\end{proof}

Theorem~\ref{Th2} allows us to use functional properties of the expected value of the switch process to 
investigate the covariance of the stationary switch process. In particular, combining Theorem~\ref{Th2} with Theorem~\ref{Th1} yields a partial characterization of the covariance functions.  

\begin{theorem}\label{Th3}
Let $C(t)$ be a symmetric function around zero, $t\in\mathbb{R}$ such that following conditions are satisfied for all $t\in[0,\infty)$,
\begin{itemize}
    \item[$(i)$] $C(t)\geqslant 0$, 
    \item[$(ii)$] $C'(t)\leqslant 0$,
    \item[$(iii)$] $C''(t)\geqslant 0$,
    \item[$(iv)$] $C(0)=1$.
\end{itemize}
Then $C(t)$ is the covariance function of a stationary switch process with a $GD(2)$ switching time distribution. 
\end{theorem}
\begin{proof}
The conditions follows directly from Proposition~\ref{Efunk}, Theorem~\ref{Th2} and the condition $C(0)=1$.
\end{proof}

If the switching time distribution also belongs to $GD(2)$, then from Corollary~\ref{Corr1} there is an explicit relation between the expected value of the switch process and the switching time distribution. A similar result exists for the stationary switch process. 

\begin{corollary}\label{Corr3}
Let $C(t)$ be the covariance of the stationary switch process and the switching time distribution belong to $GD(2)$, with the divisor $\tilde{F}$ then  
\begin{align*}
   1+\frac{\mu}{2} C'(t)&=\tilde{F}(t) \\
   \frac{\mu}{2}C''(t)&=\tilde{f}(t)
\end{align*}
\end{corollary}

Even if the switching time distribution does not belong to $GD(2)$ the expected value of the switching time distribution can still be obtained using Theorem~\ref{Th2}.  

\begin{corollary}\label{Corr4}
Let $\mu$ be the expected value of the switching time distribution, $C(t)$ the covariance of the stationary switch process, and $E(t)$ the expected value function of the corresponding switch process, then 
\begin{align*}
   \mu=2 \int_0^\infty E(u)du.
\end{align*}
\end{corollary}
\begin{proof}
From Theorem~\ref{Th2}, we have $ C'(t)=-\frac{2}{\mu}E(t)$ and for some constant $\alpha$
\begin{align*}
    -\frac{2}{\mu}\int_0^t E(u)du&=C(t)+\alpha. \\
\end{align*}
We obtain $\alpha=-1$ by taking  $t\rightarrow 0^+$ and using $C(0)=1$. Then
\begin{align*}
    -\frac{2}{\mu}\int_0^\infty E(u)du&=-1.
\end{align*}
\end{proof}

\section{Examples}\label{examples}
\noindent
To illustrate the switch process and the relation between covariance, expected value, and the switching time distribution three examples are presented. The first is a very simple example used to illustrate the basic concepts, the second shows the utility of Theorem~\ref{Th2} and the third employs the conditions of Theorem~\ref{Th3} to derive the switching time distribution. 

\subsection*{Exponential switching time} Consider a switch process with exponential switching times with intensity $\lambda$. It is well known that this distribution belongs to $GD(\infty)$ and therefore also belongs to $GD(2)$. We derive the expected value function through the Laplace transform using Proposition~\ref{LEXt}. In particular, 
\begin{align*}
\mathcal{L}(E)(s)=\frac{1}{s}\frac{1-\frac{\lambda}{\lambda+s}}{1+\frac{\lambda}{\lambda+s}}=\frac{1}{2\lambda+s},
\end{align*}
which corresponds to 
\begin{align*}
E(t)=e^{-2\lambda t}.
\end{align*}
We note that $E(t)$ is non-negative and decreasing for all $t\geqslant 0$, which is in agreement with the results of Theorem~\ref{Th1}. Let us now derive the covariance function of the corresponding stationary switch process using Theorem~\ref{Th2}. Clearly $\mu=\frac{1}{\lambda}$ and we have 
\begin{align*}
    C'(t)=-\frac{2}{\frac{1}{\lambda}}e^{-2\lambda t}=-2\lambda e^{-2\lambda t}. 
\end{align*}
By solving for $C(t)$ and computing $C''(t)$ it is clear that the conditions of Theorem~\ref{Th3} are satisfied. 

\subsection*{Gamma switching time} 
\begin{figure}[H] \label{fig1}
\includegraphics[width=0.95 \textwidth]{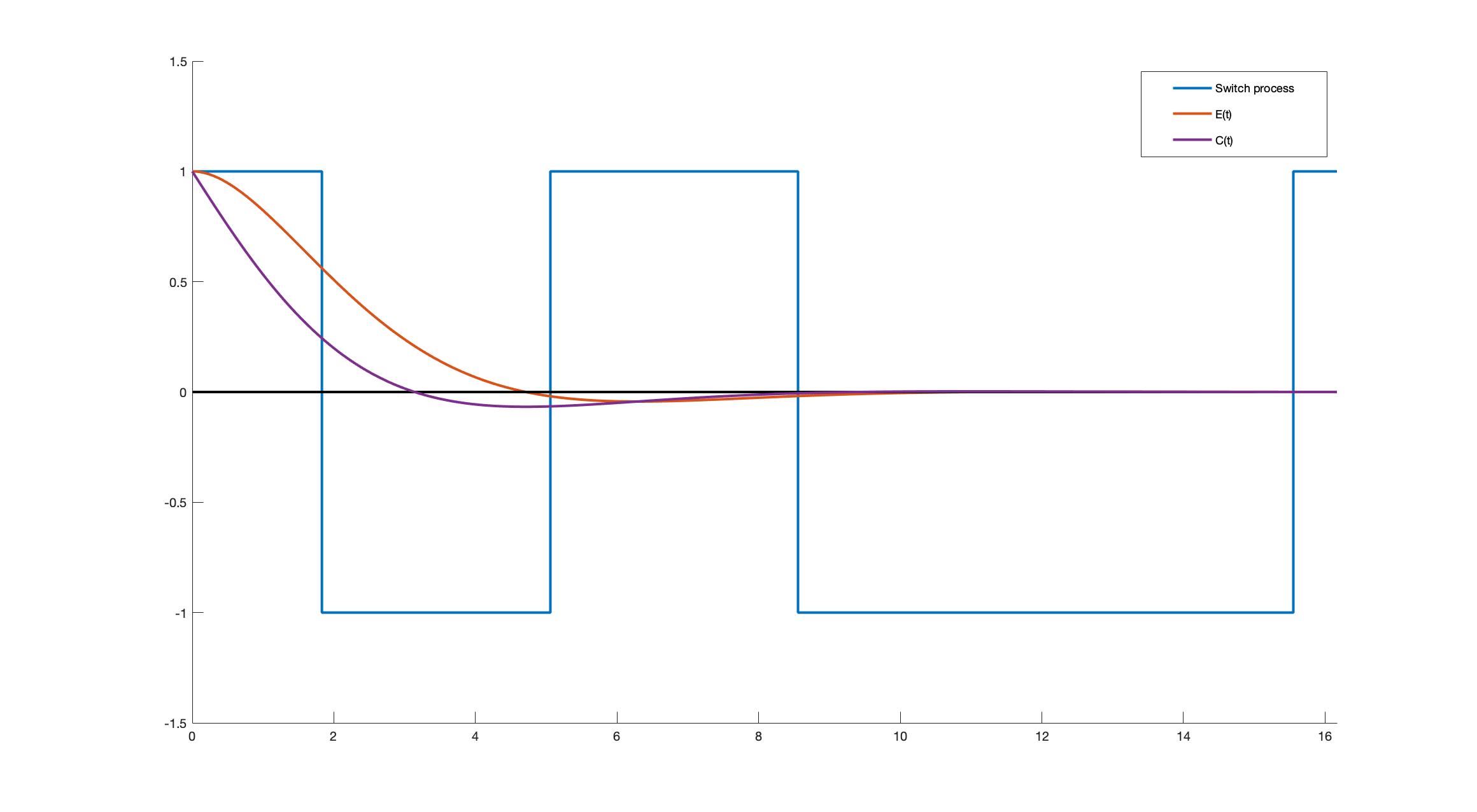}
\caption{Sample path of switch process generated from $T \in \Gamma(2,2)$, $E(t)$ and $C(t)$ of its stationary counterpart.}
\label{figgamma}
\end{figure}

Consider instead a process with gamma distributed switching times, with parameters $\theta=2$ and $k=2$. This distribution does not  belong to the class of $GD(2)$ distributions, which can be verified through Proposition~\ref{LGD}. Computing the expected value function, through the Laplace transform, we have 
\begin{align*}
    \mathcal{L}(E(t))(s)&=\frac{1}{s} \frac{1-(1+2s)^{-2}}{1+(1+2s)^{-2}}=\frac{2+2s}{1+2s+2s^2}, \\
    E(t)&=\sqrt{2} \sin \left(\frac{2t+\pi }{4} \right) e^{-\frac{t}{2}}. 
\end{align*}
Since $E(t)$ is oscillating, the condition $(ii)$ of Theorem~\ref{Th1} is not satisfied. This is not surprising since \cite{ThinG} has shown that a count process with gamma-distributed arrival times, $k>1$ cannot be obtained as the thinning of some other count process. In the context of this paper, the divisor cannot be obtained and therefore the connection between geometric divisibility and the functional properties of the expected value cannot explicitly be derived in a straightforward manner from Theorem \ref{Th1}. 
However, the covariance function of its stationary counterpart can still be obtained through Theorem~\ref{Th2}. Since the switching time distribution is known, $\mu$ is in this case equal to four. Solving the differential equation with the condition $C(0)=1$, we have for $t \geqslant 0$
\begin{align*}
    C'(t)&=- \frac{1}{\sqrt{2}}\sin\left(\frac{2t+\pi }{4}\right) e^{-\frac{t}{2}}, \\
    C(t)&=\cos\left( \frac{t}{2}\right)e^{-\frac{t}{2}}. 
\end{align*}
This example illustrates the utility of Theorem~\ref{Th2}, even if the switching time distribution does not belong to $GD(2)$. In Figure~\ref{fig1} a realization of the switch process with gamma switching times and its corresponding $E(t)$ and $C(t)$, is shown.

\subsection*{GD(2) switching time from covariance} In the previous two examples, we started with a switching time distribution and derived $E(t)$ of the switch process and $C(t)$ of the corresponding stationary switch process. Let us instead investigate if the function 
\begin{align*}
    h(t)=\frac{2}{\pi} \arcsin \left( \frac{1}{\cosh(\frac{t}{2})} \right),
\end{align*}
is a valid covariance function of a stationary switch process. We start by computing the first and second derivative of $h(t)$:
\begin{align*}
    h'(t)&=-\frac{1}{\pi}\frac{\tanh(\frac{t}{2}) \rm{sech}(\frac{t}{2})}{\sqrt{1-{\rm sech}^2(\frac{t}{2})}}=-\frac{1}{\pi}\rm{sech}\left(\frac{t}{2}\right), \\
    h''(t)&=\frac{1}{2 \pi} \tanh\left( \frac{t}{2}\right) \rm{sech}\left( \frac{t}{2}\right). 
\end{align*}
From this, the conditions of Theorem~\ref{Th3} are satisfied and the function $h(t)$ corresponds to the covariance of a specific stationary switch process, with a $GD(2)$ switching time distribution. From Corollary~\ref{Corr4} $\mu$ can be obtained and combining this with Corollary~\ref{Corr3} the switching time divisor is obtained 
\begin{align*}
   \tilde{F}(t)&=1-\rm{sech}\left( \frac{t}{2}\right), \\
    \tilde{f}(t)&= \frac{1}{2} \tanh\left( \frac{t}{2}\right) \rm{sech}\left( \frac{t}{2}\right).
\end{align*}
To obtain the explicit switching time distribution is difficult using the standard method, through the Laplace transform. Numerical methods can be used to obtain the inverse of the Laplace transform. An alternative approach is to simulate from the divisor by using the inverse sampling method and the stochastic representation of $2$-geometric divisible distributions to simulate from the full distribution. 

\section{Applications} \label{appli} 
\noindent
The  concept of persistency is of interest in statistical physics. Persistency is related to the exceedance time distribution of a stochastic process, i.e, the time between an up crossing of some level and a down crossing of the same level. These, in turn, lead to Palm distributions of the point process of instants of the level crossings.  To obtain the explicit exceedance time distribution, or the Palm distribution,  is a well-known difficult problem in probability theory, see \cite{LittlewoodOfford,Rice,https://doi.org/10.1111/sjos.12248}. Therefore a \textit{persistency coefficient}, which describes the tail behavior of the exceedance time distribution is often sought after instead. However, even obtaining exact and analytical persistency coefficients is still to a large extent an open problem. For example, the persistency coefficient was explicitly found for diffusion processes in dimension two by \cite{Exact2018} and then only for the zero level crossings. Therefore approximation methods are commonly used for the purpose.  One of them is the independent interval approximation \textit{(IIA)}. For an extensive overview see for example \cite{IIA}. 

In essence, the IIA approach works as follows. For the stochastic process of interest, define the clipped process by computing the sign of the process. The time intervals when the process is one or minus one will not be independent. However, for simplicity, one can assume that the dependency is negligible and treat them as independent. This process can be viewed as a stationary switch process and the covariance can then be used to infer information from the switching time distribution. The covariance of the clipped process is directly related to the covariance of the process of interest and has an explicit form if the underlying process is Gaussian. The covariance function of the clipped process is then matched to the covariance of the stationary switch process. This is then used to approximate the persistency coefficient. 

Two fundamental questions related to this approach should be posed. The first is if the approximated distribution is indeed a valid probability distribution, i.e. the IIA is mathematically sound. The second is if an explicit form of this distribution can be obtained. 

The covariance function of a clipped Gaussian process is 
\begin{align*}
    C(t)=\frac{2}{\pi} \arcsin\left(r(t)\right),
\end{align*}
where $r(t)$ is the covariance function of the process which is clipped, i.e our process of interest. 
By computing $C''(t)$, conditions on $r(t)$ can be obtained for when $C(t)$ corresponds to a $GD(2)$ distribution. We have the subsequent proposition which follows directly from Theorem~\ref{Th3}.
\newpage
\begin{proposition}\label{propapli}
Let $r(t)$ be the covariance function of a zero-mean Gaussian process such that the following conditions are satisfied for $t\geqslant0$
\begin{itemize}
    \item[(i)] $r(t)\geqslant0,$
    \item[(ii)] $r'(t)\leqslant0,$
    \item[(iii)] $r''(t)\geqslant-(r'(t))^2 r(t) / (1-r(t)^2).$ 
\end{itemize}
Then the approximated switching time distribution, using the IIA method, belong to the class of $GD(2)$ distributions. 
\end{proposition}

This proposition answers the posed questions and provides a partial answer to when the IIA method provides a valid probability distribution. The conditions of Proposition~\ref{propapli} are easy to verify and the conditions are satisfied by a large set of functions. 

Consider the diffusion process in two dimensions. The covariance of the clipped process has the following form 
\begin{align*}
    C(t)=\frac{2}{\pi} \arcsin \left( \frac{1}{\cosh(\frac{t}{2})} \right).
\end{align*}
From the last example in Section~\ref{examples} this corresponds to a $GD(2)$ switching time distribution. Since the distribution of the divisor is known, it is possible to infer more information than only the persistency coefficient, such as the full approximated distribution and its properties. It might require some numerical methods but requires fewer computations than, for example simulating trajectories and then estimating the exceedance time from these simulations.

\section{Conclusions} 
\noindent
To characterize which functions correspond to the expected value of the switch process is a difficult problem. By exploring the relationship between the functional properties of the expected value and the class of $2$-geometric divisible distributions, a partial answer to the problem is given.  

An explicit relation between the expected value function of the switch process and the covariance function of the stationary switch process is presented. It leads to corresponding relations between the $2$-geometric divisible switching time distributions and the covariance of the stationary switch process. It enables the recovery of the switching time distribution from the covariance function under some conditions which are easy to verify. This constitutes a partial solution to the well-known problem of obtaining the switching time distribution from the covariance function of a continuous-time binary process. 

The complete answers to both the above-mentioned problems are still unknown. However, the partial answers given in this work provide an explicit answer for an important class of functions and switching time distribution. 

\section*{Acknowledgment} \noindent 
The author is very grateful for many fruitful discussions with Krzysztof Podgórski, helpful inputs on geometric divisibility from Tomasz Kozubowski, and the inspiration and direction from  Georg Lindgren. The author is also thankful for the helpful comments and support from Yvette Baurne and Joel Danielsson.
Financial support of the Swedish Research Council (VR) Grant DNR: 2020-05168 is acknowledged.

\newpage

\appendix
\section{Proofs} \label{proofs} \noindent
In this appendix, the proofs of Section~\ref{preliminaries} are collected in order of appearance. Some of these results may be previously known, however, to keep the paper self-contained proofs are provided. 
Some supporting results are needed, but first a remark on notation. Let $f(t)$ be the  probability density function of a random variable with support on the positive real line, the Laplace transform of its corresponding probability measure is defined by $\Psi_F(s)=\int_0^\infty e^{-ts}dF(t)$. The Laplace transform of the cumulative distribution function of the same random variable will be denoted by $\mathcal{L}(F)(s)=\int_{0}^\infty  e^{-ts}F(t)dt$. The two are related through $\mathcal{L}(F)(s)=\Psi(s)_F/s$, since $F$ is the convolution of the density function with the Heaviside step function.  We use $\ast$ for convolutions of functions and $\star$ is used to denote the convolution of probability distributions i.e. the distribution function of the sum of random variables.  

\begin{lemma}\label{fgsup}
Let $G$,$H$ be distribution function on $[0,\infty)$ with corresponding  densities $g,h$. Then for each $t\geqslant 0$ 
\begin{align*}
    (g \ast h)(t) \leqslant \sup_{u>0}g(u)H(t).
\end{align*}
\end{lemma}
\begin{proof}[Proof of Lemma~\ref{fgsup}] 
Since $g,h$ are densities on $[0,\infty)$, $g(t-x)=0$ for $x>t$, we have
\begin{align*}
    (g\ast h)(t)&=\int_{-\infty}^\infty g(t-x)h(x)dx=\int_{0}^t g(t-x)h(x)dx \\
    &\leqslant  \int_{0}^t \sup_{u>0}g(u) h(x)dx=\sup_{u>0}g(u)H(t).
\end{align*}
\end{proof}

\begin{corollary}\label{corrbound}
Let $G$ be a distribution function on $[0,\infty)$, with the corresponding density $g$. Then $\forall n,m \in \mathbb{N}$ we have 
\begin{align*}
    g^{\ast (n+m)}(t) \leqslant \sup_{u>0}g(u)^{\ast n} G(t)^{m\star}. 
\end{align*}
\end{corollary}

\begin{lemma}\label{HGbound}
Let $G$ and $H$ be a distribution function of $[0,\infty)$, with corresponding densities $g,h$, then for $t\geqslant0$
\begin{align*}
    (H\star G)(t) \leqslant G(t) H(t).
\end{align*}
\end{lemma}
\begin{proof}[Proof of Lemma~\ref{HGbound}] 
The result follows from 
\begin{align*}
     (H\star G)(t)&=\int_{-\infty}^\infty \int_{-\infty}^{t-x_1} g(x_1) h(x_2)dx_1dx_2 \\
     &=\int_{-\infty}^\infty g(x_1) H(t-x_1)dx_1=\int_0^t g(x_1)H(t-x_1)dx_1 \\
     &\leqslant H(t)\int_0^t g(x_1)dx_1=H(t)G(t).
\end{align*}
\end{proof}
\noindent
From the associative property of convolutions and Lemma~\ref{HGbound}, the following bound holds
\begin{align*}
    H^{n \star }\leqslant H^n(t).
\end{align*}
The proofs of the section now follow in the order they appeared. 

\begin{proof}[Proof of Proposition~\ref{EXt}] 
From the law of total expectation and since $P(N(t)=n)=F^{n \star } - F^{ (n+1) \star }$ we have,
\begin{align*}
EX(t)&=E \left[ E[X(t) \vert N(t)] \right] =\sum_{n=0}^\infty (-1)^n \left( F^{n \star } - F^{(n+1) \star } \right)(t). 
\end{align*}
Which can also be expressed as
\begin{align*}
E(t)&= \sum_{n=0}^\infty (-1)^n F^{n \star } + \sum_{n=0}^\infty (-1)^{n+1}F^{(n+1) \star }(t)=1+2\sum_{k=0}^\infty (-1)^k F^{k \star}(t).
\end{align*}
\end{proof}

\begin{proof}[Proof of proposition~\ref{LEXt}]
We omit technicalities but from Proposition~\ref{EXt} and the relation $\mathcal{L}(F)(s)=\Psi_F(s)/s$, we have for $s>0$
\begin{align*}
\mathcal{L}(E)(s)&= \sum_{n=0}^\infty (-1)^n \mathcal{L}(F^{n \star } - F^{ (n+1) \star})(s) \\
&=\sum_{n=0}^\infty (-1)^n \frac{\Psi_F^n(s)-\Psi_F^{n+1}(s)}{s} \\
&= \frac{1}{s} ( 1- \Psi_F(s)) \sum_{n=0}^\infty (-\Psi_F(s))^n \\
&=\frac{1}{s} \frac{ 1- \Psi_F(s)}{ 1+\Psi_F(s)}.
\end{align*}
The second equation of the proposition is obtained by solving the equation above for $\Psi_F(s)$.
\end{proof}

\begin{proof}[Proof of Proposition~\ref{limitXt}]
From the last expression of the proof of Proposition~\ref{EXt}, we have
\begin{align*}
E(t)=1+2 \sum_{n=1}^\infty (-1)^n F^{n \star}(t).
\end{align*}
It is enough to show that the infinite series converges to zero when $t$ goes to zero. Using Lemma~\ref{HGbound} we obtain 
\begin{align*}
\bigg \vert \sum_{n=1}^\infty (-1)^n F^{n \star}(t) \bigg \vert \leqslant \sum_{n=1}^\infty F^{n}(t)= \frac{F(t)}{1-F(t)}, ~t\ge 0.
\end{align*}
Since $\lim_{t \rightarrow 0} F(t)=0$, the limit  at zero is obtained. 
\\
The limit at infinity follows from the Key Renewal Theorem (see page 86 in \cite{theoryofpoint1}) as follows. 

We first note that 
\begin{align*}
    E(t)=2P(X(t)=1)-1
\end{align*}
and thus it is sufficient to show that
\begin{align*}
\lim_{t\rightarrow \infty} P(X(t)=1)=\frac{1}{2}.
\end{align*}
The line of the argument goes along a standard technique in the renewal theory. 
We start with the renewal equation that is easy to show by conditioning on the first two switch intervals that
\begin{align*}
P(X(t)=1)&=1-F(t)+\int_0^t P(X(t-u)=1) ~dF^{2\star}(u).
\end{align*}

It follows from Lemma~4.1.I (Renewal Equation Solution), p.~69 in  \cite{theoryofpoint1} that 
\begin{align*}
P(X(t)=1)=1-F(t)+\sum_{k=1}^\infty \int_0^t 1- F(t-u) ~dF^{k\star}(u).
\end{align*}
Then by the Key Renewal Theorem, see p.86 in  \cite{theoryofpoint1}
\begin{align*}
\lim_{t\rightarrow \infty} P(X(t)=1)=
\frac{1}{2\mu}\int_0^\infty 1-F^{2\star}(u)~du=\frac 1 2.
\end{align*}
where $\mu$ is the expected value of the switching time distribution.
\end{proof}

\begin{proof}[Proof of Proposition~\ref{EXtprim}] 
Consider the partial sum
\begin{align*}
E_N(t)=1+2\sum_{k=1}^N (-1)^k F^{k \star}(t),
\end{align*}
which, from the last equation of the proof of Proposition~\ref{EXt}, is well defined for each $N$ and converges to $E(t)$ when $N\rightarrow \infty$. From Theorem 7.17, \cite{RudinPoP}, it is sufficient to show local uniform convergence of $E'_N(t)$, in order to exchange summation and derivation. Thus for each $0< t_0 < t_1$, we need to show that for each $\epsilon >0 $ there exist a $N$ such that 
\begin{align*}
    \big\vert E'_N(t)-E'(t) \big\vert  < \epsilon
\end{align*}
for all $t\in [t_0,t_1]$.

For given $t_0$ and $t_1$, let $t\in [t_0,t_1]$, and choose an $N>l$.
Then for every $\epsilon>0$, we have  
\begin{align*}
     \big\vert 2 \sum_{k=1}^N (-1)^k f^{\ast k}(t) - 2 \sum_{k=1}^\infty (-1)^k f^{\ast k}(t) \big\vert &= \big\vert 2 \sum_{k=N+1}^\infty (-1)^k f^{\ast k}(t) \big\vert \\
     = \big\vert 2 \sum_{n=1}^\infty (-1)^{N+n} f^{\ast l + (N-l)+n }(t) \big\vert & 
     \leqslant 2 \sum_{n=1}^\infty f^{\ast l}*f^{\ast (N-l)+n }(t).
\end{align*}
Let us first assume that the support of $f$ is $(0,\infty)$. By applying Corollary \ref{corrbound} and Lemma \ref{HGbound} we bound the last term in the above by 
\begin{align*}
     2\sup_{u>0} f^{\ast l}(u) \sum_{n=1}^\infty  F^{*(N-l+n)}(t)
     & \le
     2\sup_{u>0} f^{\ast l}(u)
     \sum_{n=1}^\infty F^{N-l+n}(t_1)
     \\
    &= 
    2\frac{F(t_1)}{1-F(t_1)} \sup_{u>0} f^{\ast l}(u) F^{N-l}(t_1), 
\end{align*}
where $F(t_1)<1$ since $f(t)$ is supported on $(0,\infty)$.
For each $\epsilon>0$, there will always exist a $N$ such that 
\begin{align*}
    F^{N-l}(t_1) <  \epsilon \ \frac{1-F(t_1)}{2 F(t_1) \sup_{u>0} f^{\ast l}(u)},
\end{align*}
which proves uniform convergence of $E'_N(t)$ on compact subsets of $\mathbb{R}^+$.

The proof extends to the case of a distribution with support in the positive real line by considering how the support expands for each convolution and by grouping the terms in such a way that $t_1$ falls within the support of the iterated convolutions, i.e.  $F^{k \star}(t_1)<1$, for some $k>0$, where $k$ is then used as the size of grouping the subsequent terms in the series. 
\end{proof}

\begin{proof}[Proof of Proposition~\ref{LGD}]
If $F\in GD(r)$  it has the following Laplace transform, described in Section~\ref{preliminaries},
\begin{align*}
    \Psi_F(s)&=\frac{\frac{1}{r}\Psi_{\tilde{F}}(s)}{1-(1-\frac{1}{r})\Psi_{\tilde{F}}(s)}.
\end{align*}
Solving for $\Psi_{\tilde{F}}(s)$, we have
\begin{align*}
    \Psi_{\tilde{F}}(s)=\frac{r\Psi_F(s)}{1 +(r-1)\Psi_F(s)}. 
\end{align*} 
If and only if $F\in GD(r)$, $\Psi_{\tilde{F}}(s)$ will be a completely monotone function. This follows directly from Bernstein's theorem (see Theorem 1, page 415 in \cite{FellerV2}). 
\end{proof}

\begin{proof}[Proof of Proposition~\ref{setlemma}]
Let $X\in GD(r)$ with the distribution function $F$ and corresponding divisor $\tilde{F}$, then we have 
\begin{align*}
    \Psi_F=\frac{\frac{1}{r}\Psi_{\tilde{F}}(s)}{1-(1-\frac{1}{r})\Psi_{\tilde{F}}(s)}.
\end{align*}
Substituting the above into the expression of Proposition~\ref{LGD}, we have
\begin{align*}
    \frac{u\Psi_F(s)}{1 +(u-1)\Psi_F(s)}
    &=\frac{u\frac{\frac{1}{r}\Psi_{\tilde{F}}(s)}{1-(1-\frac{1}{r})\Psi_{\tilde{F}}(s)}}{1 +(u-1)\frac{\frac{1}{r}\Psi_{\tilde{F}}(s)}{1-(1-\frac{1}{r})\Psi_{\tilde{F}}(s)}}\\
    &=\frac{\frac{u}{r} \Psi_{\tilde{F}}(s) }{ 1-(1-\frac{1}{r})\Psi_{\tilde{F}} (s)+(u-1)\frac{1}{r}\Psi_{\tilde{F}} (s)} \\
    &=\frac{\frac{u}{r}\Psi_{\tilde{F}}(s)}{1-(1-\frac{u}{r})\Psi_{\tilde{F}}(s)}.
\end{align*}
If $\frac{u}{r}<1$ or, equivalently, $u<r$, then the last expression is completely monotone since it corresponds to $\frac{u}{r}-$geometric summation of the divisor of $X$. From Proposition~\ref{LGD}, $X$ is therefore both $r$-geometric and $u-$geometric divisible. 
\end{proof}

\newpage
\bibliographystyle{apalike}
\bibliography{References.bib}




 
\end{document}